\documentclass[12pt]{article}
\input epsf
\usepackage{epsfig}
\usepackage{amsmath}

\usepackage{xcolor}

\begin{document}
\def\R{\mathbf{R}}
\def\bC{\mathbf{\overline{C}}}
\def\C{\mathbf{C}}
\def\l{\ell}
\def\intr{\mathrm{int\ }}
\title{Electrostatic skeletons}
\author{Alexandre Eremenko\thanks{Supported by NSF grant.},$\,$
 Erik Lundberg,
 and Koushik Ramachandran\thanks{Supported by NSF grant DMS-1162070.}}
\maketitle
\begin{center}
{\em Dedicated to the memory of Andrei Gonchar and Herbert Stahl}
\end{center}
\begin{abstract}
Let $u$ be the equilibrium potential of a
compact set $K$ in $\R^n$.
An electrostatic skeleton of $K$
is a positive measure $\mu$
such that the closed support $S$ of $\mu$ has connected complement and no interior,
and the Newtonian (or logarithmic, when $n=2$) potential of $\mu$
is equal to $u$ near infinity.
We prove the existence and uniqueness of an
electrostatic skeleton for any simplex.

{\em MSC: 31A12, 31A25. Key words: potential, equilibrium, subharmonic
function, inverse problem, analytic continuation.}
\end{abstract}

Let $K \subset \R^n$ be a compact set which is
regular for the Dirichlet problem.
The equilibrium potential $u_K$ of $K$ with pole at infinity
is a positive harmonic function in $\R^n \setminus K$ which satisfies
\begin{equation}\label{norm}
u_K(x)=1 + O(1/|x|^{n-2}),\quad x\to\infty,\quad\text{for}\quad n>2,
\end{equation}
or
$$u_K(z)=\log|z|+O(1),\quad z\to\infty,\quad\text{for}\quad n=2,$$
and whose boundary values on $K$ are zero. Thus for $n=2$ the equilibrium
potential is the same as the Green function of $\R^2\backslash K$
with pole at infinity. We use the unusual normalization in (\ref{norm})
to have a closer analogy between the cases $n=2$ and $n>2$:
our equilibrium potentials, when extended
by $0$ on $K$, become subharmonic functions
in $\R^n$.

A positive measure $\mu$ with closed support $S \subset K$ is called
an {\em electrostatic skeleton} of $K$ if $S$ has empty interior,
$\R^n\backslash S$ is connected,
and
$$u_K(x)= 1 - \int \frac{1}{|x-y|^{n-2}} d\mu(y), \quad x \in\R^n \setminus K, \quad \text{for } n>2,$$
or, if $n=2,$
$$u_K(z)=\int\log|z-\zeta|d\mu + \gamma, \quad z\in\R^2\backslash K, $$

for some constant $\gamma.$

\vspace{0.1in}

For example, an ellipsoid has an electrostatic skeleton
supported on its $(n-1)$-dimensional focal ellipsoid \cite{LK}.

For $n>2$, if
$$P(x) = \sum_{i=k}^{d} \frac{a_k}{|x-x_k|^{n-2}},\quad a_k>0,$$
then the discrete measure with atoms of mass $a_k$
at each point $x_k$ is an electrostatic skeleton for the level sets
$\{ x:P(x)\geq L\}$, $L>0$.
The analogous example for $n=2$
is that polynomial lemniscates
have skeletons supported on finite sets.

An electrostatic skeleton of a special region bounded by two circles
in $\R^2$
appears in \cite{LSS}.
It is clear that there are sets $K$ which do not have an electrostatic
skeleton. For example, a Jordan curve in $\R^2$ whose boundary is nowhere
analytic. If $u$ is the equilibrium potential of such a $K$, then
any level set $\{ z:u(z)<c\}$ also does not have an electrostatic skeleton.

Motivated by his study of the asymptotic behavior of
zeros of Bergman (area orthogonal) polynomials, E. B. Saff proposed the problem about the
existence of the skeletons in 2003 and mentioned it at several conferences
(see also the reference to Saff in the recent paper by Lundberg
and Totik \cite{LT}).
In particular, Saff asked whether every convex polygon has an electrostatic skeleton.

Here we prove the existence and uniqueness of skeletons
for triangles (and splices in higher dimensions).

An analogous question with potentials of
the volume or surface area measure instead of
the equilibrium potential was considered by B. Gustafsson~\cite{G}.
\vspace{.1in}

\noindent
{\bf Theorem 1.} {\em If $K$ is a simplex, then there exists a unique electrostatic skeleton.}
\vspace{.1in}

{\em Proof.} Let $L_j,\; j=1,\ldots,n+1$ be the open $(n-1)$-faces of $K$, and $\l_j$ the
reflection in $L_j$. Then the equilibrium potential$$u_K(z)=\int\log|z-\zeta|d\mu+c, \quad z\in\R^2\backslash K, \quad \text{for } n=2,$$
-where $c$ is a constant.
$u=u_K$ extends by reflection:
$$u_j=-u\circ\l_j,\quad j=1,\ldots,n+1.$$
As $K$ is convex, these functions $u_j$ are negative harmonic functions
in the interior of $K$. The boundary values of $u_j$ on $\partial K$
are zero on the closure
$\overline{L_j}$ and strictly negative on  $K\backslash\overline{L_j}$.
It follows that
the $u_j$ are pairwise distinct negative harmonic functions in $\intr K$. So there is no open set $V\subset K$ where $u_j(x)=u_k(x)$ for $x\in V$
and $j\neq k$.

\vspace{0.1in}

Next, define a function $w$ by 
$$w(x)=\left\{\begin{array}{ll} u(x),& x\in\R^n\backslash K,\\
\max\{ u_1(x),\ldots,u_j(x)\},& x\in K.\end{array}\right.$$ 

Clearly, $w$ is subharmonic in $\R^n$ and, harmonic in $\R^n\backslash K.$  We claim that the skeleton is the Riesz measure $\mu,$ of $w.$ It is evident that the potential of $\mu$ matches the equilibrium potential outside $K.$ 
So it remains to prove that $\R^n\backslash S[w]$ is connected, where $S[w]$ is the support of the Riesz measure $\mu.$

\vspace{0.1in}

\noindent The support can explicitly described as  
$$ S[w] = \{x\in K: u_i (x) = u_j(x),\hspace{0.05in}\mbox{for some}\hspace{0.05in} i\neq j\}.$$

Suppose that $\R^n\backslash S[w]$ has a bounded component $D$,
then $D\subset K.$ To every component $D$ of $\intr K\backslash S[w]$ corresponds a
number $k(D)$ such that $w(x)=u_k(x),\; x\in D$. Suppose without loss of generality that
$k(D)=1$, that is $u_1(x)>u_j(x),\; z\in D,\; j\in\{2,3,..,n+1\}$. Let $a$ be a point in $D$.
Let $R$ be the straight line segment connecting a point $b\in L_1$
with the vertex opposite to $L_1$ which passes through $a$.
Let $R_2, R_3,..,R_{n+1}$ be reflections of $R$ with respect to the other faces
$L_2,L_3,..,L_{n+1}$. These segments all lie outside $K$.
Introduce the direction on $R,R_2,R_3,..,R_{n+1}$ from
their common vertex. See Fig.~1.
We claim that our original equilibrium potential
$u$
is strictly increasing on these segments.
\begin{figure}
\begin{center}
\includegraphics[height=6cm]{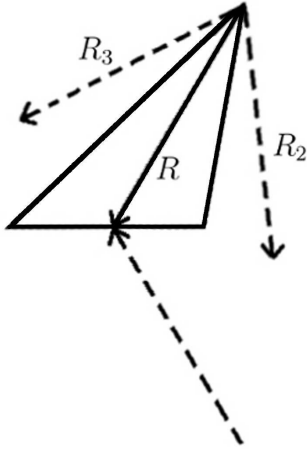}
\caption{Construction for the proof of Theorem 1.}
\label{figG}
\end{center}
\end{figure}

\vspace{0.1in}

This follows from the fact that level hypersurfaces
of the equilibrium potential
of a convex set are convex (for $n=2$, see \cite{P}, and for $n>2$,
\cite{Lewis}), and a straight line
can cross a convex hypersurface at most twice; if the straight line
segment begins inside a convex hypersurface and ends outside, then
it crosses the boundary only once.
Therefore, our segments
$R_2,R_3,..,R_{n+1}$ cross each level surface at most once, so the function
$u$ is increasing on $R_2,R_3,..,R_{n+1}$. Therefore, $u_j=-u\circ\l_j$
with $j\in\{2,3,..,n+1\}$ are decreasing on $R$, while $u_1$ is increasing
on $R$. So as we have $u_1(x)>\max\{u_2(x),u_3(x),..,u_{n+1}(x)\}$ for
$x=a$ and $x=b$ we conclude that this inequality holds on the
whole segment $[a,b]\subset R$ giving a contradiction, which proves that
$S_1$ does not divide $\R^n.$ This proves the existence of a skeleton for simples. Uniqueness follows from  Proposition $1$ below.

\vspace{.1in}

\noindent
{\bf Remarks and  conjectures.}
\vspace{.1in}

\noindent
1. It is easy to see that a non-convex (Jordan) polygon in $\R^2$ cannot
have a skeleton.
\vspace{.1in}

\noindent
2. Now we give an example of a circular triangle which does not
have an electrostatic skeleton. Let $U$ be the unit disc, and $D\subset
\Delta=\bC\backslash\overline{U}$ the hyperbolic triangle with
vertices at $1,\exp(\pm2\pi i/3)$. Then repeated reflections of $D$
in the sides tile $\Delta$. A more familiar picture of this tiling
is obtained by changing the variable to $1/z$. Therefore the equilibrium
potential has an analytic continuation to $\Delta\backslash E$, where $E$
is a discrete set of logarithmic singularities. As these singularities are
dense on the unit circle, $K=\C\backslash D$ cannot have an
electrostatic skeleton.

\vspace{.1in}

\noindent 3. For any convex polytope, define $S[w]$ as in Theorem $1.$ We conjecture that $S[w]$ does not divide $\R^n.$ 
This implies existence of a skeleton as in the Proof of Theorem $1,$ and it also implies uniqueness (Proposition $1$ below).

\vspace{0.1in}

In general, the electrostatic skeleton of a Jordan region $K$ does not have
to be unique. An example of non-uniqueness was
constructed by Zidarov \cite{Z}, it is also reproduced in \cite{G}.
For simplicity we only discuss this example for $n=2$.
Zidarov constructed two distinct probability measures whose supports
are trees consisting of finitely many straight segments, and whose potentials
coincide in a neighborhood of $\infty$. These potentials are
equal in a neighborhood of $\infty$ with the potential of
the area measure of a non-convex polygon.
The level sets of these potentials
$\{ z:u(z)\leq c\}$ with sufficiently large $c$ are Jordan regions
with non-unique skeletons.
\begin{figure}[h]
\begin{center}
\includegraphics[height=4cm]{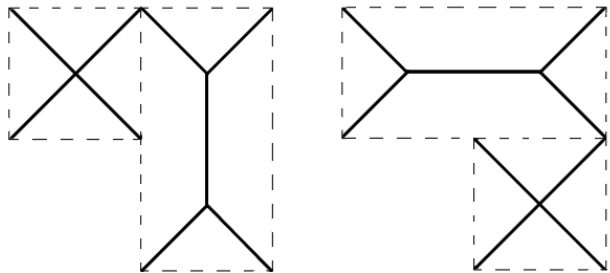}
\caption{Two trees supporting positive measures whose potentials coincide near infinity.}
\end{center}
\end{figure}
\vspace{.1in}

\noindent

4. When the equilibrium potential in the definition of skeleton
is replaced by the volume measure on $K$, the resulting notion is
called a ``mother body''. Gustafsson \cite{G} proved existence and
uniqueness of the mother body for every convex polytope in $\R^n$.

\vspace{.1in}

\noindent
{\bf Proposition 1.} {\em Let $K$ be a convex polytope in $\R^n.$ Suppose that $K$ admits an electrostatic skeleton. Let $w$ be defined as in Theorem $1.$ Using the notation from the proof of Theorem 1, if $S[w]$ does not divide $\R^n$, then
the electrostatic
skeleton is unique, and its closed support is $S[w]$.}
\vspace{.1in}

{\em Proof.}
By assumption, electrostatic skeleton exists. Let $S$ be its support, and $v$ the potential of the
skeleton. Clearly $S\subset K$. By assumption, $S[w]$ does not divide
$\R^n$. Then the complement
$\intr K\backslash S[w]$ consists of $N$ components, $D_j$, where $N$ is the
number of $(n-1)$-faces of $\partial K$, such that $D_j$
contains $L_j$.

As $S$ is closed and has empty interior, the set $G=\intr K\backslash S$
is open and non-empty. Let $p$ be a point in
$G$. As $S$ does not divide space, there is a curve
$\gamma$ in $\R^n\backslash S$ starting from $p$ and ending outside $K$,
and $v$ has an analytic continuation on this curve. Consider the
point $q$ where $\gamma$ leaves $K$ for the first time, and let $q\in L_k$
(it is clear that $q$ belong to a face of dimension $n-2$ of $K$).
As $u_k$ is the immediate analytic continuation of $u$ to
the interior of $K$, we conclude that $v(x)=u_k(x)$ in a neighborhood of $p$.
Therefore
\begin{equation}\label{W}
v(x)\in\{ u_1(x),\ldots, u_N(x)\}\quad\mbox{for every}\quad x\in K,
\end{equation}

In particular, $v$ is continuous.
Recall that $u_k(x)>u_j(x)$ for $x\in D_k$ and all $j\neq k$. As $v$ is
continuous, it follows from (\ref{W}) that $v(x)=u_k(x)$
for $x \in D_k$. So $v=w$
by continuity, and $S=S[w]$.

\vspace{.1in}

We thank Edward Saff and Vilmos Totik for an interesting discussion
of this paper.

\vspace{0.1in}

{\em Department of Mathematics

Purdue University

West Lafayette IN 47907-2067

\vspace{0.1in}

eremenko{@}math.purdue.edu

elundber{@}math.purdue.edu

kramacha{@}math.purdue.edu}

\end{document}